\documentclass{amsart}
\usepackage{amsmath}
\usepackage{graphicx}
\usepackage{amsfonts}
\usepackage{amssymb}
\usepackage{hyperref}
\newtheorem{theorem}{Theorem}
\theoremstyle{plain}

\newtheorem{lemma}{Lemma}

\newtheorem{proposition}{Proposition}

\numberwithin{equation}{section}

\begin{document}
\title{One Dimensional Conformal Metric Flows}
\author  { Yilong Ni and
Meijun Zhu}
\address{Department of Mathematics\\
The University of Oklahoma\\
Norman, OK 73019\\
}

\begin{abstract}
In this paper we continue our studies of  the one dimensional
conformal metric flows, which were introduced in \cite{NZ1}. We
prove the global existence and convergence of the flows, as well
as the exponential convergence of the metrics under these flows.
\end{abstract}

\maketitle

\section{Introduction}
In the past two decades, there have been intensive studies of the
evolutions of plane curves in ${\bf R}^2$. Such research has broad
impacts in the development of mathematics as well as computer
vision and human beings' recognition system.  The simplest model
seems to be the evolution of a plane curve along its normal
direction in ${\bf R}^2$:
\begin{equation}
{\bf x}_t=\phi(k){\bf n},
\label{0-1}
\end{equation}
where ${\bf n}$ is the inner unit norm of the curve ${\bf x}$, $k$
is its curvature and $\phi(\cdot)$ is a given function. For
$\phi(x)=x$, it is the well-known {\it curve shortening flow}, see
for example, Gage \cite{G} and  Gage and Hamilton \cite{GH}; For
$\phi(x)=x^{1/3}$, it is equivalent to the affine flow, which was
studies by Sapiro and Tannenbaum \cite{ST}, and by Alvarez, et al.
\cite{Al}. In the past decade, the affine flow has its great
successful applications in image processing, mainly because of its
larger (comparing with curve shortening flow) invariant group.

In \cite{NZ1} we propose a novel approach to the evolution of a
given convex planar curve. This approach seems more intrinsic and
natural from conformal geometric point of view. In fact, we will
show in this paper, among other things, that it yields a simple
proof of the global existence of the affine flow.

We recall from \cite{NZ1} that if $(S^1,g_s)$ is the unit circle
with the induced metric $g_s=d \theta \otimes d \theta$ from ${\bf
R}^2$, for any metric $g$ on $S^1$ (for example, this metric could
be given by reparametrizing the circle), we write $g:=d \sigma
\otimes d \sigma=v^{-4}g_s$ for some positive function $v$,  and
then define a general {\it $\alpha$- scalar curvature} of $g$ for
any positive constant $\alpha$  by
$$R^\alpha_g=v^3(\alpha v_{\theta\theta}+v).$$
Thus $ R^\alpha_{g_s}=1. $ The corresponding {\it $\alpha$-
conformal Lapalace-Beltrami} operator of $g$ is defined by
$$
L^\alpha_g=\alpha\Delta_g+R^\alpha_g,
$$
where $\Delta_g=D_{\sigma \sigma}$. It is proved in \cite{NZ1} that
$L^\alpha_g$ is a conformal covariant. More precise, the following
proposition is proved:
\begin{proposition}
For $\varphi>0$, if $g_2=\varphi^{-4}g_1$ then
$\kappa_{g_2}=\varphi^3 L^\alpha_{g_1}\varphi$, and
$$
L^\alpha_{g_2}(\psi)=\varphi^3 L^\alpha_{g_1}(\psi\varphi),
\quad\forall \psi\in C^2({\bf S^1}).
$$
\label{prop0-1}
\end{proposition}
The general   {\it $\alpha$-scalar curvature} flow can be defined
as
\begin{equation}
\partial_t g=-R_g^\alpha g.
\label{0-2}
\end{equation}
It was pointed out in \cite{NZ1} that  two cases of $\alpha=1$ and
$\alpha=4$ are of special interest: First of all,
 the affine curvature of any given convex curve can be
viewed as a {\it $1$-scalar curvature}. To see this, let ${\bf
x}(\theta)\subset{\bf R}^2$ ($\theta\in [0, 2\pi]$) be a convex
simple closed curve parameterized by the angle $\theta$ between
the tangent line and $x$-axis. One can introduce a new parameter
of the affine arc-length  $\sigma$ by requiring
$$[{\bf x}_\sigma, {\bf x}_{\sigma \sigma}]=1.
$$
This can be done by defining
$$
\sigma(\theta)=\int_0^{\theta} k^{-2/3} d \theta=\int_0^{s}
k^{1/3} d s,
$$
where $s$ is the parameter of the arc length, $k=k(\theta)$ is the
curvature of the curve. Let $v=k^{1/3}$, we have
$$
g_1:=d\sigma \otimes d\sigma=v^{-4}d\theta \otimes d\theta=v^{-4}g_s.
$$
Then the affine curvature of ${\bf x}(\theta)$ is given by:
$$
\kappa=v^3(v_{\theta\theta}+v),
$$
which is the same as $R^1_g$. On the other hand, given $(S^1,g)$ with
$g=u^{-4}g_s$, we may define ${\bf C}_u(\theta)\subset {\bf R}^2$ as
$$
{\bf C}_u(\theta)=\left(\int_0^\theta u^{-3}(\theta)\cos\theta d\theta,
\int_0^\theta u^{-3}(\theta)\sin\theta d\theta \right).
$$
If $u(\theta)$ satisfies the orthogonal condition
\begin{equation}
\label{on} \int_0^{2\pi}u^{-3}(\theta)\cos\theta d\theta
=\int_0^{2\pi}u^{-3}(\theta)\sin\theta d\theta=0,
\end{equation}
then ${\bf C}_u(\theta)\subset {\bf R}^2$ so defined is a convex simple
closed curve and its affine curvature is equal to $R^1_g$.

The affine flow, which is successfully used in  image processing
(see \cite{ST}), is defined as
$$
{\bf x}_t(\sigma,t)={\bf x}_{\sigma\sigma}.
$$
It is known (see, for example, \cite{ST}) that the above flow is
equivalent to
$$
\partial_t g=-\kappa_g g.
$$
We will see in the following shortly that the affine flow is
equivalent to the one-dimensional normalized affine flow:
\begin{equation}\label{affine}
\partial_t g=(\overline{\kappa}_g-\kappa_g) g,\,\, L(0)=2\pi.
\end{equation}
 where $ \overline{\kappa}_g={\int \kappa_g d\sigma}/{\int
d\sigma}$ and $L(t)=\int d\sigma$.

For the case of $\alpha=4$, the {\it $4$-scalar curvature} $R^4_g$ can
be view as an analog one-dimensional Yamabe flow. Let $g=d \delta
\otimes d \delta=u^{-4} g_s =u^{-4} d \theta \otimes d \theta$ for
some positive function $u \in C^1(S^1)$. We define the scalar
curvature $k_g$ of metric $g$ by
$$
k_g:=u^3(4u_{\theta \theta}+u),$$ and the corresponding  {\it conformal
Lapalace Beltrami} operator by
\begin{equation} A_g:=4\Delta_g
+k_{g},  \label{0-3}
\end{equation}
that is $A_g=L^4_g$. Thus $A_g$ is a conformal covariant by
Proposition \ref{prop0-1} and
$$
k_{g_2}=\varphi^3 L^\alpha_{g_1}\varphi,
$$
for  $g_2=\varphi^{-4}g_1$. The one-dimensional  Yamabe flow is guided by
\begin{equation}
\partial_t g=-k_g g.
\label{0-4}
\end{equation}
Again this flow is equivalent to the normalized
Yamabe flow:
\begin{equation}
\partial_t g=(\overline{k}_g-k_g) g, \,\, L(t)=2\pi,
\label{0-5}
\end{equation}
where $\overline{ k}_g= {\int k_g d\sigma}/{\int d\sigma}$.

The steady states of these flows were studied in details in
\cite{NZ1}. Based on these results,  we shall prove
\begin{theorem}
For an abstract curve  $({S^1}, {u_0}^{-4}g_s)$, there is a unique
smooth solution $g(t)$ to the flow equation (\ref{0-5}) for $t\in
[0, +\infty). $ Moreover, $g(t) \to g_\infty$ in $L^\infty(S^1)$
exponentially fast in $t$ as $t \to +\infty$, and the 4-scalar
curvature of $g_\infty$ is constant. \label{theorem1}
\end{theorem}

And for the affine flow, we shall prove

\begin{theorem}
Let  $(S^1, g_0)$ be an abstract curve. If
$g_0=u^{-4}(\theta,0)g_s$ and $u(\theta,0)$ satisfies the
orthogonal condition (\ref{on}). Then there is a unique smooth
solution $g(t)$ to the flow equation (\ref{affine}) for $t\in
[0,+\infty)$. Moreover, $g(t) \to g_\infty$ in $L^\infty(S^1)$
exponentially fast in $t$  as $t \to +\infty$, and the 1-curvature
of $g_\infty$ is constant. \label{theorem2}
\end{theorem}

In the case of curve shortening flow, the exponential convergence
of the curvatures (thus the  curves in  ${\bf R}^2$) is proved in
\cite{GH}. For the affine flow, the convergence of the affine
curvatures is obtained in \cite{ST} (see also,  \cite {A1} and
\cite{A2}). Our approach is quite different. In fact,
as being mentioned in \cite{NZ1}, we are motivated by the study of
conformal geometry and  the early work of X. Chen \cite{C1}, where
he initiates the study of Calabi flow via integral estimates. In
the same spirit we expect to pursue our systematical study of the
analog $\alpha$-scalar curvature problems on high dimensional
spheres in forthcoming papers.

We organize the paper as follows. In Section 2, we derive some
basic properties about the flows and prove the global existence of
the flows. We then deal with the Yamabe flow in Section 3. We
first derive the $L^\infty$ bound for the metric, then we prove
the exponential convergence of the curvatures and the metrics.  We
devote the whole Section 4 to the proof of the exponential
convergence of metrics for 1-scalar curvature flow. For readers'
convenience, we list the needed sharp inequalities from \cite{NZ1}
in the appendix.

\medskip
\section{Basic properties and global existence}
We shall first show that the general flow equation (\ref{0-2}) for
$g(t)$ is equivalent to a normalized flow.  In fact, if we choose
$$\hat g(t)=\frac{4\pi^2}{L^2(0)}\exp\left(\int_0^t\overline{R}_g^\alpha(\tau)
d\tau\right) g(t),
$$
where $L(t)=\int d \sigma(t)$ and $g(t)=d \sigma(t) \otimes d
\sigma(t)$, and  a new time variable
$$\hat t=\frac{4\pi^2}{L^2(0)}\int_0^t \exp\left(\int_0^\delta
 \overline{R}_g^\alpha(\tau)d \tau \right)d \delta,
$$
then equation (\ref{0-2}) can be written as
$$
\partial_{\hat t} \hat g= (\frac{\int \hat R^\alpha_{\hat g}d \hat \sigma}{\int
d\hat  \sigma}-\hat R^\alpha_{\hat g}) \hat g.$$

From now on, we shall focus on  the normalized flow:
\begin{equation}
\partial_t g=(\overline{ R}^\alpha_g -R^\alpha_g)g, \,\, L(0)=2\pi.
\label{1-1}
\end{equation}
\begin{lemma}
Along the conformal flow equation (\ref{1-1}) with
$g(\sigma,t)=u^{-4}g_s$, ${R}:=R^\alpha_g$ satisfies
\begin{equation}
{R}_t= \frac{\alpha}4 \Delta {R}+ {R}({R} -\overline{R}), \label{1-2}
\end{equation}
the metric satisfies
\begin{equation}\label{1-3}
\partial_t (d \sigma)=\frac 12 (\overline{R}-R) d \sigma,
\end{equation}
and  $u$ satisfies
\begin{equation}
u_t=\frac 14(R -\overline{R})u.  \label{1-4}
\end{equation}
 \label{lem1-1}
\end{lemma}

\begin{proof}
$$
(\overline{R}-R) g=\partial_t g=(-4)u^{-5}u_tg_s =(-4)u^{-1}u_tg,
$$
that is, $u_t=\frac 14( R -\overline{R})u$. Thus
$$
\partial_t (d\sigma)=(u^{-2}d\theta)_t=-2u^{-3}u_td\theta
=\frac12(\overline{R}-R )d\sigma.
$$
Using the conformal invariance of $L^\alpha_g$, we have
\begin{align*}
R_t=&(u^3L^\alpha_{g_s}u)_t=3u^2u_tL^\alpha_{g_s}u+u^3L^\alpha_{g_s}(u_t)
=3u^{-1}u_t R+u^3L^\alpha_{g_s}\left(u \frac {u_t}u\right)\\
=&\frac34(R -\overline{R})R+L^\alpha_{g}(u_t/u)
=\frac {\alpha }4\Delta_g R+R(R -\overline{R}).
\end{align*}
\end{proof}

It follows from (\ref{1-3}) that
$$
\partial_t\int_{S^1} d\sigma=\int_{S^1}\frac12(\overline{R^\alpha_g}
-R^\alpha_g)d\sigma=0.
$$
Thus flow (\ref{1-1}) preserves the  arc length with respect to
metric $g$ (i.e. $\int_0^{2\pi} u^{-2} d \theta=L(0)=2\pi$).
Moreover, along the flow, we see from the following lemma that the total
curvature is strictly increasing unless $R^\alpha_g$ is a constant.

\begin{lemma}Along flow (\ref{1-1}), we have
\begin{equation}
\partial_t \overline{R}^\alpha_g= \frac 1{4\pi}
\int_{S^1}(R^\alpha_g-\overline{R}^\alpha_g)^2 d\sigma. \label{1-5}
\end{equation}
 \label{lem1-2}
\end{lemma}
\begin{proof}
\begin{align*}
\partial_t\overline{R}^\alpha_g=&\frac1{2\pi}\int_{S^1} (R^\alpha_g)_td\sigma
+\frac1{2\pi}\int_{S^1} R^\alpha_g\partial_t (d\sigma)\\
=&\frac1{2\pi}\int_{S^1} R^\alpha_g (R^\alpha_g-\overline{R}^\alpha_g)d \sigma
+\frac1{4\pi}\int_{S^1} R^\alpha_g(\overline{R}^\alpha_g-R^\alpha_g)d \sigma\\
=&\frac1{4\pi}\int_{S^1} R^\alpha_g(R^\alpha_g-\overline{R}^\alpha_g)d \sigma
=\frac1{4\pi}\int_{S^1}( R^\alpha_g-\overline{R}^\alpha_g)^2d\sigma\ge0.
\end{align*}
\end{proof}

We note that (\ref{1-4}) can also be written as a heat equation
\begin{equation}
u_t= \frac{\alpha}4 u^4 \Delta_{g_s} u+ \frac {1}4 u^5 -\frac 14
\overline{R}^\alpha_g u. \label{1-6}
\end{equation}

We are now ready to prove the global existence for $\alpha=4$
(which we refer to the {\it Yamabe flow}) and $\alpha=1$ (which we
refer to the {\it affine flow}). For convenience we use $R_g$ to
replace $R^4_g$ and $\kappa_g$ to replace $R^1_g$  in the rest of
this paper. So the normalized Yamabe flow is written as
\begin{equation}
\partial_t g=(\overline{ R}_g -R_g)g, \ g(0)=u(\theta,
0)^{-4} g_s; \label{Yequ}
\end{equation}
And the normalized affine flow is written as
\begin{equation}
\partial_t g=(\overline{\kappa}_g -\kappa_g)g,  \ g(0)=u(\theta,
0)^{-4} g_s. \label{Aequ}
\end{equation}

We first show the global existence of the Yamabe flow via the
following local $L^\infty$ estimate.

\begin{proposition} Suppose $g(t)=u^{-4}g_s$ satisfies (\ref{Yequ}).
Then for any given $t_0>0$, there is a positive
constant $c=c(t_0)>0$ such that $$ \frac 1{c(t_0)} \le u( t)\le
c(t_0),\qquad t\in[0,t_0].$$ \label{prop1-1}
\end{proposition}
\begin{proof}
For simplicity, we use $R$ to replace $R_g$ in this proof.  From
(\ref{1-2}) we know that
$$ R_t+\overline{R} R \ge \Delta R.$$
It follows from the maximum principle that \begin{equation} R
\ge \min_\sigma R (\sigma,0) \cdot e^{-\int_0^t
\overline{R } d \tau}.\label{1-7}
\end{equation}
Due to Theorem 2 in \cite{NZ1} (Theorem B in the appendix) we also
know that $\overline{R }\le 1$. Thus there is a  constant
$c_1(R(\sigma,0))$, such that
\begin{equation} R(\sigma, t)
\ge c_1(R(\sigma,0)),\qquad t\in[0,t_0].\label{1-8}
\end{equation}
It then follows from (\ref{1-4}) that  we have
\begin{equation}
u(\sigma, t)=u(\sigma, 0) \cdot e^{\frac 14
\int_0^t(R-\overline{R}) d \tau} \ge c_2(R(\sigma,0),
t_0)>0,\quad t\in[0,t_0]. \label{1-9}
\end{equation}

To estimate the upper bound for $u(\sigma, t)$, we first observe
that $u$ satisfies
$$
4 u_{\theta\theta}+ u=R u^{-3},\quad u>0, \ \mbox{  and } \int_0^{2\pi}
u^{-2}d\theta=2\pi.
$$
Multiplying the above by $u$ and then integrating it over $S^1$,
we obtain
\begin{equation}\label{1-10}
\int_0^{2\pi} u^2d\theta-4 \int_0^{2\pi}u_\theta^2d\theta
=\int_0^{2\pi} R u^{-2}d\theta=2 \pi \overline{R} \ge 2 \pi \overline{R}_0.
\end{equation}
Let $M(t)=|\{\theta\mbox{ : }u(\theta,t)\ge 2\}|$. Then (\ref{1-9})
implies
\begin{align*}
2\pi=\int_0^{2\pi} u^{-2}d\theta=&\int_{u\ge 2}u^{-2}d\theta
+\int_{u<2}u^{-2}d\theta \\
\le&\frac{M(t)}4+(2\pi-M(t))c_2(R(\sigma,0),t_0)^{-2}.
\end{align*}
Therefore there exists $\delta(t_0)>0$, such that $2\pi-M(t)\ge
\delta(t_0)$. That is
$$
|\{\theta\mbox{ : }u(\theta,t)\le 2\}|>\delta(t_0),\quad
\mbox{for} \  t\in[0,t_0].
$$
If  $\sup_{t \in [0, t_0) }\int_0^{2\pi} u^2(t)d\theta =\infty$,  then
there exists a sequence $ t_i \to t_*\le t_0$, such that $\int_0^{2\pi}
u^2(t_i)d\theta=\tau_i^2\to\infty$  as $i\to \infty$. We define
$v_i=u(t_i)/|\tau_i|$. It follows from (\ref{1-10}) that $v_i$
satisfies
$$
\int_0^{2\pi} v^2_id\theta=1,\mbox{ and } 4 \int_0^{2\pi}(v_i)_\theta^2
d\theta\le\int_0^{2\pi}  v_i^2 d\theta -\frac{2 \pi
\overline{R}_0}{\tau_i^2}\le c_3,
$$
which yields that $\{v_i\}$ is a bounded set in
$H^1\hookrightarrow C^{0,\frac12}$. Therefore up to a subsequence
$v_i \rightharpoonup v_0$ in $H^1$ weekly and $v_0 \in
C^{0,\frac12}$. It is easy to see that $v_0$ satisfies
$|\{\theta\mbox{ : }v_0(\theta)=0\}| \ge\delta(t_0)$, and
$$
0\le  \int_0^{2\pi} v_0^2d\theta-4
\int_0^{2\pi}(v_0)_\theta^2d\theta.
$$
 Without loss of
generality, we can assume that the north pole of $S^1$ is in
$\{\theta\mbox{ : }v_0(\theta)=0\}$. Let $v^*$ be the
symmetrization of $v_0$, such that $v^*(\mbox{north pole})=0$. Let
$\phi: S^1\to \mathbf R$ be the stereographic projection and
$w(y)=\frac2{1+y^2}\cdot v^*\circ\phi^{-1}(y)$. Since
$|\{\theta\mbox{ : }v_0(\theta)=0\}|\ge\delta(t_0)$, it is easy to
see that there exist $R>0$ such that $w(y)=0$ for $|y|>R$ and
$$
0 \le\int_0^{2\pi} v_0^2d\theta-4 \int_0^{2\pi}(v_0)_\theta^2d\theta\\
=-4\int_{-R}^R|w'(y)|^2dy,
$$
contradiction. Therefore $\int u^2d\theta$ is bounded on
$[0,t_0]$, so is $\int u_\theta^2d\theta$. Thus $u(t)$ is
bounded in $H^1\hookrightarrow C^{0,\frac12}$, which implies that
 there exists a $c(t_0)>0$ such that
$$
\frac{1}{c(t_0)}\le u(t) \le c(t_0),\qquad t\in[0,t_0].
$$
\end{proof}

Since $u(t)$ satisfies (\ref{1-6}), we know from the standard
parabolic estimates that $u(t)$ exists for all $t\in[0,+\infty)$,
this completes the proof of the global existence for (\ref{Yequ}).

\medskip

For the affine flow, we need to modify the proof slightly. We
first show that the orthogonality is preserved under the affine
flow.

\begin{lemma}
 Suppose that $g(t)=u^{-4}g_s$ satisfies (\ref{Aequ}). If
 $\int_0^{2\pi}\cos \theta \cdot u^{-3}(\theta,0)d \theta= \int_0^{2\pi}\sin \theta \cdot u^{-3}(\theta,0)d
 \theta=0,$ then for all $t>0$,
 \begin{equation}
\int_0^{2\pi}\cos \theta \cdot u^{-3}(\theta,t)d \theta=
\int_0^{2\pi}\sin \theta \cdot u^{-3}(\theta,t)d
 \theta=0.\label{1-13}
 \end{equation}
\label{lem1-3}
\end{lemma}

\begin{proof}
From (\ref{1-4}) and the definition of $\kappa$ we have
\begin{align*}
\partial_t \int_0^{2\pi}\cos \theta \cdot u^{-3}(\theta,t)d
\theta&=-3 \int_0^{2\pi}\cos \theta \cdot u^{-4}(\theta,t)
u_t(\theta,t)d \theta\\
&=-\frac 34 \int_0^{2\pi} \kappa \cos \theta \cdot
u^{-3}(\theta,t)d \theta+\frac{3 \overline{\kappa}}4
\int_0^{2\pi}\cos \theta \cdot
u^{-3}(\theta,t)d \theta\\
&=\frac{3 \overline{\kappa}}4 \int_0^{2\pi}\cos \theta \cdot
u^{-3}(\theta,t)d \theta.
\end{align*}
Thus
$$
 \int_0^{2\pi}\cos \theta \cdot u^{-3}(\theta,t)d
\theta=C e^{\int_0^t \overline{\kappa}(\tau) d \tau}.
$$
Since  $\int_0^{2\pi}\cos \theta \cdot u^{-3}(\theta,0)d
\theta=0$, we have $C=0$, thus $ \int_0^{2\pi}\cos \theta \cdot
u^{-3}(\theta,t)d \theta=0$. Similarly, we can obtain $
\int_0^{2\pi}\sin \theta \cdot u^{-3}(\theta,t)d \theta=0$.

\end{proof}

We are ready to obtain the local $L^\infty$ estimate for the
metric under the affine flow.
\begin{proposition}
Suppose that $g(t)=u^{-4}g_s$ satisfies (\ref{Aequ}) and the
initial data satisfies
 $\int_0^{2\pi}\cos \theta \cdot u^{-3}(\theta,0)d \theta= \int_0^{2\pi}\sin \theta \cdot u^{-3}(\theta,0)d
 \theta=0.$
 Then for any given $t_0>0$,
there is a positive constant $\tilde c=\title c(t_0)>0$ such that
$$ \frac 1{\tilde c(t_0)} \le u( t)\le \tilde c(t_0),\qquad
t\in[0,t_0].$$ \label{prop1-2}
\end{proposition}
\begin{proof} For simplicity, we use $\kappa$ to replace $\kappa_g$ in the proof.
In the same spirit to the proof of Proposition \ref{prop1-1}, we
first derive  the low bound for $\kappa$.

From (\ref{1-2}) we know that
$$ \kappa_t+\overline{\kappa} \kappa \ge \frac{1}4 \Delta \kappa.$$
It follows from the maximum principle that $$ \kappa \ge
\min_\sigma \kappa (\sigma,0) \cdot e^{-\int_0^t \overline{\kappa
} d \tau}.$$ Due to Theorem 1 in \cite{NZ1}(Theorem A in the
appendix) we also know that $\overline{\kappa}\le 1$. Thus there
is a constant $\tilde c_1(\kappa(\sigma,0))$, such that
\begin{equation} \kappa(\sigma, t)
\ge \tilde c_1(\kappa(\sigma,0)),\qquad t\in[0,t_0].\label{1-14}
\end{equation}
From (\ref{1-4}) (also using (\ref{1-14})) we have
\begin{equation}
u(\sigma, t)=u(\sigma, 0) \cdot e^{\frac 14
\int_0^t(\kappa-\overline{\kappa}) d \tau} \ge \tilde
c_2(\kappa(\sigma,0), t_0)>0,\quad t\in[0,t_0]. \label{1-15}
\end{equation}

To estimate the upper bound for $u(\sigma, t)$, we first note
that $u$ satisfies
$$
u_{\theta\theta}+ u=\kappa u^{-3},\quad u>0 \mbox{  and  }\int_0^{2\pi}
u^{-2}d\theta=1.
$$
Multiplying $u$ and integrating over $S^1$, we obtain
\begin{equation}\label{1-16}
\int_0^{2\pi} u^2d\theta- \int_0^{2\pi}u_\theta^2d\theta
=\int_0^{2\pi} \kappa
u^{-2}d\theta=2 \pi\overline{\kappa} \ge 2 \pi\overline{\kappa}_0,
\end{equation}
where $\overline{\kappa}_0=\overline{\kappa}(\sigma,0)$. Let
$\tilde M(t)=|\{\theta\mbox{ : }u(\theta,t)\ge 2\}|$. Then
(\ref{1-15}) implies
\begin{align*}
2\pi=\int_0^{2\pi} u^{-2}d\theta=&\int_{u\ge 2}u^{-2}d\theta
+\int_{u<2}u^{-2}d\theta \\
\le&\frac{\tilde M(t)}4+(2\pi-\tilde M(t)) \tilde
c_2(\kappa(\sigma,0),t_0)^{-2}.
\end{align*}

Therefore there exists $\tilde \delta(t_0)>0$, such that
$2\pi-\tilde M(t)\ge \tilde \delta(t_0)$. That is
$$
|\{\theta\mbox{ : }u(\theta,t)\le 2\}|>\tilde \delta(t_0),\quad
\mbox{for} \  t\in[0,t_0].
$$
If  $\sup_{t \in [0, t_0) }\int_0^{2\pi} u^2(t)d\theta =\infty$,  then
there exists a sequence $t_i \to T_* \le t_0$, such that $\int_0^{2\pi}
u^2(t_i)d\theta=\tau_i^2\to\infty$  as $i\to \infty$. We define
$v_i=u(t_i)/|\tau_i|$. It follows from (\ref{1-16}) that $v_i$
satisfies
$$
\int_0^{2\pi} v^2_id\theta=1,\mbox{ and }  \int_0^{2\pi}(v_i)_\theta^2
d\theta\le\int_0^{2\pi}  v_i^2 d\theta -\frac{2 \pi
\overline{\kappa}_0}{\tau_i^2}\le c_3,
$$
which means that $\{v_i\}$ is a bounded set in $H^1\hookrightarrow
C^{0,\frac12}$. Therefore up to a subsequence $v_i \rightharpoonup
v_0$ in
 $H^1$ weekly and $v_0 \in C^{0,\frac12}$. It is easy to see that $v_0$ satisfies
$|\{\theta\mbox{ : }v_0(\theta)=0\}| \ge \tilde \delta(t_0)$, and
\begin{equation}
0\le  \int_0^{2\pi} v_0^2d\theta-
\int_0^{2\pi}(v_0)_\theta^2d\theta. \label{1-17} \end{equation}
Moreover, since $u(\theta, 0)$ satisfies (\ref{1-13}), $u(\theta,
t_i)$ and $v_i(\theta)$ must also satisfy (\ref{1-13}) by Lemma
\ref{lem1-3}. It follows that if $[a, b]\subset [0, 2 \pi)$ is a
set such that
\begin{equation}
b>a, \ v_0(a)=v_0(b)=0, \ v(\theta)>0 \mbox{ for }\theta \in (a, b),
\label{1-18} \end{equation}
 then $b-a<\pi$ (see more details in the proof of Theorem 1 in \cite{NZ1}). Let
$$\{\theta\mbox{ : }v_0(\theta)>0\}=\cup_k I_k,$$
where $I_k=(a_k, b_k)$ are all intervals satisfying (\ref{1-18}).
For each interval $(a_k, b_k)$, we define
$$w_k(\theta):=
\left\{ \begin{array}{rll} & v_0( \frac\theta2+a_k) \ \ \ &
\mbox{for} \ \ \theta\in [0,2(b_k-a_k)] \\
 & 0 \ \ \ & \mbox{for \ other }\ \theta. \end{array} \right.
$$
It is easy to see that $w_k \in H^1(S^1)$, and
$$
\int_0^{2 \pi} (w_k^2- 4(w_k)_\theta^2)d\theta=
\int_{a_k}^{b_k} (v_0^2- (v_0)_\theta^2)d\theta.
$$
Similar argument to that in the proof of Proposition
\ref{prop1-1}, we have for all $k$,
$$
\int_0^{2 \pi} (w_k^2- 4(w_k)_\theta^2)d\theta >0.
$$
Thus
$$\int_0^{2 \pi} (v_0^2-(v_0)_\theta^2)d\theta=
\sum_{k}\int_{a_k}^{b_k} (v_0^2- (v_0)_\theta^2)d\theta >0.
$$
Contradiction to (\ref{1-17}).
 Therefore $\int u^2d\theta$ is bounded on
$[0,t_0]$, so is $\int u_\theta^2d\theta$. Thus $u(t)$ is
bounded in $H^1\hookrightarrow C^{0,\frac12}$, which implies that
there exists a constant $c(t_0)>0$ such that
$$
\frac{1}{c(t_0)}\le u(t) \le c(t_0),\qquad t\in[0,t_0].
$$
\end{proof}

Again, since $u(t)$ satisfies (\ref{1-6}), we know from the
standard parabolic estimates that the affine flow $u(t)$ exists
for all $t\in[0,+\infty)$.

\medskip

\section{convergence for the Yamabe flow}

In this section, we shall prove the exponential convergence for
the metric under Yamabe flow. We will achieve this goal by two
steps. First, we shall use moving plane method to establish
 the uniform bound for metrics; Then we will obtain the exponential convergence via
 energy estimates and a Kazdan-Warner type identity.

 We recall: along the Yamabe flow
(\ref{Yequ}), $u$ satisfies
\begin{equation}\label{3-1}
u_t=u^4(u_{\theta\theta}+\frac u4) -\frac14\overline{R}_g u,\mbox{
on }S^1\times[0,\infty),
\end{equation}
and $R_g$ satisfies \begin{equation}\label{3-2} R_t=\Delta
R+R(R-\overline{R}), \end{equation} where and throughout this
section, we use $R$ to replace $R_g$ and $\Delta$ to replace
$\Delta_g$.

Along the same line in Ye's classical paper \cite{Ye}, we shall
first use moving plane method to prove the following Harnack
inequality.
\begin{proposition}
There exists constant $C>0$, such that
$$
\frac1C\le\frac{|\nabla_{g_s}u|}{u}\le C, \quad \forall
t\in[0,\infty).
$$
\label{prop3-1}
\end{proposition}

\begin{proof} Let $\Phi$: ${\bf x}=(x_1, x_2) \in S^1 \to y \in R^1$ be the  stereographic
projection given by
$$
x_1=\frac {2 y}{1+y^2},
$$
and
$$
x_{2} =\frac{y^2-1}{y^2+1}.
$$
Around the north pole $p_0=(0,1)$, we choose the following
coordinates: $$G(y)=(\frac {2 y}{1+y^2}, \frac{1-y^2}{y^2+1}), \ \
\ \ \ |y|<1.$$ In fact, for $0<|y|<1$, $G(y)=\Phi^{-1}(y/|y|^2),$
and $G(0)=p_0.$

Note that
$$d \theta \otimes d \theta=\sum_{i=1}^2 dx_i^2=(\frac 2{1+y^2})^2 d y \otimes d
y:=\varphi^{-4} g_0,$$ where $\varphi(y)=\sqrt{(1+y^2)/2}.$

Let $$ w(y, t)= u(\Phi^{-1}(y), t) \cdot \varphi (y).$$ Then $w(y,
t)$ satisfies
\begin{equation}\label{mv-1}
-\frac13(w^{-3})_t=w_{yy} - \frac14\overline{R}_g w^{-3} \ \ \ \ \mbox{in} \ \
(-\infty, +\infty)\times[0,+\infty).
\end{equation}

We denote $a_0(t)=u(G(0), t)=u(p_0, t)$, $a_1=u'(G(0), t),$ and
$a_{11}(t)=u'' (G(0), t).$ Then direct computations yield the
uniformly (for any given time $T>0$) asymptotic behavior of
$w(y,t)$ near infinity on ${\bf R}$:
\begin{lemma} As $y \to \pm\infty$,
$$
w(y,t)=\frac {|y|}{\sqrt 2} \cdot \big {(} a_0(t)-\frac{2a_1(t)
y}{y^2}+(\frac {a_0(t)}2+2 a_{11}(t)) \cdot \frac 1{y^2} +O(\frac
1{|y|^3})\big { )},$$
 and
$$
\frac{ d w}{d y}(y,t)=\frac y{\sqrt 2 |y|}\cdot a_0(t)- (\frac
{a_0}{2\sqrt 2}+\sqrt{2}a_{11}(t)) \cdot \frac y{|y|^3}
+O(\frac 1{|y|^3}).$$ \label {mv-lem1}
\end{lemma}

Let $$ y_c(t)=\frac {2a_1(t)}{a_0(t)}.
$$
To complete the proof of Proposition \ref{prop3-1}, we need to
prove that $y_c(t)$ is uniformly (in $t$) bounded. We shall
achieve this via the method of moving planes.

For a given $\lambda>0$, let $y^\lambda=2 \lambda-y$ and
$w_\lambda (y, t)=w(y^\lambda, t).$  For $t=0$, we know from Lemma
\ref{mv-lem1} and the standard moving plane procedure (See, for
example, \cite{CGS})  that there is a $\lambda_0>0$ such that for
all $\lambda \ge \lambda_0$, $w(y, 0) \ge w_\lambda(y,0)$ whenever
$y \ge \lambda.$  Similarly, for any given $T>0$, there is a
uniform constant $\lambda_1(T)\ge \lambda_0$ such that for each
$\lambda\ge \lambda_1$,
\begin{equation}\label{mv-3}
w(y, t) \ge w_\lambda(y,t) \ \ \ \ \mbox{for}\ \ \ t \in [0, T] \
\mbox{and} \ \ \ y \ge \lambda.
\end{equation}

Define
$$I:=\{\lambda \ : \ \lambda \ge \lambda_0, \ \ \lambda>\max_{0
\le t \le T} y_c(t), \ \ w(y, t) \ge w_\lambda(y,t)\}.
$$
From (\ref{mv-3}), we immediately know that $I$ is not empty. Our
goal is to show that $I$ is open and closed in $(\lambda_0,
+\infty)$.

To show that $I$ is open, we first observe from Lemma
\ref{mv-lem1} that $w(y, t) \equiv w_\lambda(y,t)$ can never
happen for any $\lambda \ge \lambda_0$.

Suppose that $\lambda \in I$. It follows from the strong maximum
principle that
\begin{equation}
w(y, t) > w_\lambda(y,t) \ \ \ \ \mbox{for}\ \ \ t \in [0, T] \
\mbox{and} \ \ \ y > \lambda, \label{mv-4}
\end{equation}
and from the Hopf Lemma that
\begin{equation}
dw(y, t)/dy(\lambda) >0 \ \ \ \ \mbox{for}\ \ \ t \in [0, T].
\label{mv-5}
\end{equation}
For fixed $t$, we  shift the origin to $y_c(t)$. The expansion of
$w(y, t)$ in the new coordinates will be
\begin{equation}
w(y,t)=\frac {|y|}{\sqrt 2} \cdot \big {(}  a_0(t)+(\frac {
 a_0}2+ 2 a_{11}(t)) \cdot \frac 1{y^2} +O(\frac
1{|y|^3})\big { )},\label{mv-6}
\end{equation}
 and
\begin{equation}
\frac{ d w}{d y}(y,t)=\frac y{\sqrt 2 |y|}\cdot
a_0(t)+O(\frac 1{|y|^2}). \label{mv-7} \end{equation}
 The plane $y=\lambda$
becomes the plane $y=\lambda -y_c(t)$ in the new coordinates.
Since $\lambda \in I$, the standard moving plane procedure yields
that there is a $\epsilon(t)$ such that for $\tilde {\lambda}\in
(\lambda-\epsilon(t), \lambda+\epsilon(t)),$
$$
w(y, t) > w_{\tilde{\lambda}}(y,t).
$$
Note (\ref{mv-6}) and (\ref{mv-7}) are uniformly in time $t \in
[0, T]$, we know that there is a uniform $\epsilon>0$ such that
$(\lambda-\epsilon, \lambda+\epsilon)\subset I$, thus $I$ is open.

Next we show that $I$ is closed. Let $\lambda>\lambda_0$ be in the
closure of $I$. By continuity we know that
$$
w(y, t) \ge w_\lambda(y,t) \ \ \ \mbox{and} \ \ \lambda \ge
\max_{0\le t\le T} y_c(t).$$ We claim that $\lambda>\max_{0\le
t\le T} y_c(t)$, thus $\lambda \in I$, and we know that $I$ is
closed.

To verify our claim, we argue by contradiction. If
$\lambda=\max_{0\le t\le T} y_c(t)$, then $\lambda= y_c(t_0)$ for
some $t_0 \in [0,T].$ Now, we choose $y_c(t_0)$ as the new origin,
and consider the stereographic projection $\Phi: \ S^1 \to {\bf
R}$. Let $z(\theta, t)$ and $z_\lambda(\theta, t)$ be given by
$$
\Phi^*(w^{-3} g_{R})=z^{-3} g_{S^1}, \ \ \ \Phi^*(w_\lambda^{-3}
g_{R})=z_\lambda^{-3} g_{S^1}.$$ Then $z$ and $z_\lambda$ are
defined on $S^1_+\times [0,T]$, where $S^1_+=\{\theta \ : \
-\pi/2<\theta <\pi/2\}$. And both $z$ and $z_\lambda $ satisfy
(\ref{3-1}), $z\ge z_\lambda$ and $z(\pi/2, t)=z_\lambda(\pi/2,
t)$, $z(-\pi/2, t)=z_\lambda(-\pi/2, t)$ for all $t\in[0,T]$. It
follows from the expansion (\ref{mv-7}) that
$$
\frac {d z}{d \theta}(\frac {\pi}2, t_0)=\frac {d z_\lambda}{d
\theta}(\frac {\pi}2, t_0).
$$
Thus, due to the Hopf Lemma, we know that $z(\theta, t)\equiv
z_\lambda(\theta, t).$ But $z(\theta, 0)\equiv z_\lambda(\theta,
0)$ can not happen since $\lambda>\lambda_0$. Contradiction.

We have thus shown that $I=(\lambda_0, +\infty)$, which yields
that $u'/u \le \lambda_0$ uniformly in $t$ at the north pole of
$S^1$. Similarly, if we move the line from the negative side of
the $y$-axis, we can obtain that $u'/u \ge -\lambda_1$ for some
universal positive constant $\lambda_1$  uniformly in $t$ at the
north pole. Since we can arbitrarily choose the north pole,  we
complete the proof of Proposition \ref{prop3-1}.
\end{proof}

Since $\int u^{-2}d\sigma_0=1$, Proposition \ref{prop3-1} implies
\begin{equation}\label{3-3}
\frac 1C\le u\le C,\mbox{ for any } (\sigma, t)\in[0,2\pi]\times[0,\infty).
\end{equation}

For $p \ge 2$, we define
\begin{equation}\label{3-4}
F_p(t)=\int_0^{2\pi} |R-\overline{R}|^p d \sigma.
\end{equation}
Direct computation yields
\begin{align}\label{tfp}
\partial_t F_p= &-\frac {4(p-1)}{p} \int_0^{2\pi}|\nabla(|R-\overline
R|^{p/2})|^2 d \sigma + (p-\frac 12)\int_0^{2\pi}|R-\overline R|^{p}
(R-\overline R) d \sigma \\
&+p \overline R \int_0^{2\pi}|R-\overline R|^{p} d \sigma- \frac p{4\pi}
\int_0^{2\pi}|R-\overline R|^{p-2}  (R-\overline R) d \sigma \cdot
\int_0^{2\pi}|R-\overline R|^{2} d \sigma.\nonumber
\end{align}
Thus, for any fixed $p\ge 2$, there is a constant $C_1(p)$ such
that
$$
\partial_t F_p\le C_1(p)(F_{p+1}+F_p+F_p^{1+\frac 1p})-\frac
{4(p-1)}{p} \int_0^{2\pi}|\nabla(|R-\overline R|^{p/2})|^2 d \sigma.$$ It
then follows from (\ref{3-3}) and Sobolev inequality that for any
$k>1$, there is a constant $C_2(p, k)>0$(independent of $g(t)$)
such that
$$
\partial_t F_p\le C_1(p)(F_{p+1}+F_p+F_p^{1+\frac 1p})
-C_2(p, k)F_{kp}^{\frac 1k}.
$$
Note that
$$
F_{p+1} \le F_p^{\lambda\cdot\frac {p+1}{p}}
\cdot F_{3 p }^{(1-\lambda)\cdot\frac {p+1}{3p} },
$$
where $\lambda =(2p-1)/(2p+2)$. We have, via Young's inequality,
that for some positive constants $C_3(p)$ and $C_4(p)$,
\begin{equation}\label{3-6}
\partial_t F_p\le C_3(p)(F_{p}+F_p^\beta+F_p^{1+\frac 1p})
-C_4(p)F_{3 p}^{\frac 13},
\end{equation}
where $\beta=(2p-1)/(2p-3)>1$.

\begin{lemma}\label{lem3-1} For any $p \ge 2$,
$$
F_p(t) \to 0 \ \mbox{ as }\  t \to +\infty \ \ \mbox{ and } \ \
\int_0^\infty F_p(t)dt<\infty.
$$
\end{lemma}
\begin{proof}
We first claim that $\int_0^\infty F_p(t)^\delta dt<\infty$ for
some $p\ge 2$ and some $0<\delta\le 1$, implies $F_p(t)\to 0$ as
$t \to +\infty$, and $\int_0^\infty F_{3p}(t)^{\frac13}
dt<\infty$.

To prove the claim, we note that for any $\epsilon>0$, we may
choose a $t_\epsilon>0$ such that
\begin{equation}
\label{3-7}
F_p(t_\epsilon)<\epsilon\mbox{  and  }
\int_{t_\epsilon}^\infty F_p(t)^\delta dt<\epsilon.
\end{equation}
If $\epsilon<\epsilon_0=1/(3C_3(p)+1)$, we must have $F_p(t)\le 1$
for $t>t_\epsilon$. In fact, if this is not the case, let
$t_0>t_\epsilon$ be the first time such that $F_p(t_0)=1$. It
follows from (\ref{3-6}) that
\begin{align*}
F_p(t_0)-F_p(t_\epsilon)&\le C_3(p)\int_{t_\epsilon}^{t_0}(F_{p}+F_p^\beta
+F_p^{\frac {p+1}{p}})dt\\
&\le  C_3(p)\int_{t_\epsilon}^{t_0}(F_{p}^\delta+F_p^\delta+F_p^\delta)dt
\le 3C_3(p) \epsilon.
\end{align*}
Hence $\epsilon\ge 1/(3C_3(p)+1)$, contradiction. Now for
$\epsilon<\epsilon_0$, from (\ref{3-6}) and (\ref{3-7}) we have
for $t>t_\epsilon$,
$$
F_p(t) \le F_p(t_\epsilon)+3C_3(p) \epsilon <(3C_3(p)+1) \epsilon,
$$ which implies $F_p(t) \to 0$ as $t \to +\infty$. Furthermore, it follows
from (\ref{3-6}) that
$$
C_4(p)\int_{t_\epsilon}^\infty F_{3p}^{\frac13}\le
F_p(t_\epsilon)+3C_3(p)\epsilon< \infty,
$$
which implies $\int_0^\infty F_{3p}^{\frac13}<\infty$.

From Lemma \ref{lem1-2} and Theorem 1 in \cite{NZ1}(Theorem A in
the appendix), we know that
$$
\int_0^\infty F_2(t) dt < +\infty.
$$
Then the lemma follows from the claim by induction and H\"older's
inequality.
\end{proof}

From Lemma \ref{lem1-1} and \ref{lem1-2}, we know that
\begin{equation}
\partial_t (R-\overline R) =\Delta  (R-\overline R)+R (R-\overline R)
-\frac 1{4\pi} \int_0^{2\pi} (R-\overline R)^2 d \sigma.
\label{3-8}
\end{equation}
It then follows from (\ref{3-3}), Lemma \ref{lem3-1} and the
standard estimates for parabolic equation that
\begin{equation}
||R-\overline R||_{L^\infty} \to 0 \ \ \ \mbox{as} \ \ t \to +\infty.
\label{3-9}
\end{equation}

Next, we shall obtain exponential convergence. Note that
\begin{align}\label{3-10}
\partial_t \int_0^{2\pi}|R-\overline R|^{2} d \sigma
= &-2 \int_0^{2\pi}|\nabla(R-\overline
R)|^2 d \sigma + \frac 32\int_0^{2\pi}|R-\overline R|^{2} \cdot
(R-\overline R) d \sigma \\
&+2 \overline R \int_0^{2\pi}|R-\overline R|^{2} d \sigma. \nonumber
\end{align}
Since $\overline{R}(t)$ is increasing and bounded above by $1$, we may assume
that
$$
\lim_{t\to\infty}\overline{R}(t)=r_\infty.
$$
Then $||R-r_\infty||_{L^\infty}\to 0$ as $t\to +\infty$. It
follows from (\ref{Yequ}) and (\ref{3-3}) that there is a
subsequence $t_n \to \infty$, such that $ u(t_n) \to u_\infty$ in
$ C^\alpha(S^1)$; It then follows from the classical theorem of
Simon \cite{Si} that  $ u(t_n) \to u_\infty$ in $ C^\alpha(S^1)$
as $ t \to+\infty$ for all $t$, where the scalar curvature of
metric $g_{\infty}=u_\infty^{-4} d\theta \otimes d \theta$ is
constant $r_\infty$. Since $\int_0^{2\pi} u_\infty^{-2}d\theta=1$,
from Theorem 2 in \cite{NZ1} (Theorem B in the appendix), we
obtain that $r_\infty=1$.

Let $\gamma(\theta)=\int_0^\theta u_\infty^{-2}(\theta)d\theta$,
$\tilde{u}(\gamma,t)=u(\theta(\gamma),t)/u_\infty(\theta(\gamma),t)$
and $\tilde{R}(\gamma,t)=R(\theta(\gamma),t)$. From the conformal
covariance of $R$, we have
$$
\tilde{R}(\gamma,t)=\tilde{u}(\gamma,t)^3(4D_{\gamma\gamma}\tilde{u}(\gamma,t))
+\tilde{u}(\gamma,t) ).
$$
Furthermore,
\begin{align*}
\tilde{u}(\gamma,t)&\to 1,\mbox{ as } t\to\infty,\\
\tilde{u}_t(\gamma,t)&=\frac14\left(\tilde{R}(\gamma,t)
-\overline{\tilde{R}}(\gamma,t)\right)\tilde{u}(\gamma,t),\\
\tilde{R}_t(\gamma,t)&=D_{\gamma\gamma}\tilde{R}(\gamma,t)
+\tilde{R}(\gamma,t)\left(\tilde{R}(\gamma,t)-\overline{\tilde{R}}(\gamma,t)
\right),
\end{align*}
where $\overline{\tilde{R}}(\gamma,t)=\frac1{2\pi}\int_0^{2\pi}
\tilde{R}(\gamma,t)\tilde{u}(\gamma,t)^{-2}d\gamma
=\overline R (\theta(\gamma),t)$.

Therefore, without loss of generality (up to a changing of
variable), we can assume that $u_\infty=1$, that is
$g_\infty=g_s=d \theta \otimes d \theta.$ To obtain the
exponential convergence, we shall follow Chen \cite{C1} (see also
Struwe \cite{S} for more detailed computations). Let
$\{\phi_i(\theta)\}_{i=0}^{+\infty}$ be a $L^2$-orthonormal basis
of eigenfunctions for $-\Delta_s$ on $(S^1, g_s)$ with eigenvalue
$\mu_i$ for $i=0, 1, 2, 3, ...$. We can expand
$R-C=\sum_{i=1}^\infty R^i \phi_i$. Thus
$$
-\Delta_s (R-\overline{R})=\sum_{i=1}^\infty \mu_i R^i \phi_i.$$
It follows from (\ref{tfp}) that
\begin{align} \label{3-12}
\partial_tF_2=&-2\int_0^{2\pi}|\nabla(R-\overline R)|^2 d \sigma+\frac 32
\int_0^{2\pi}|R-\overline R|^{2} (R-\overline R) d \sigma  +2 \overline R
\int_0^{2\pi}|R-\overline R|^{2} d \sigma\\
=&-2\int_0^{2\pi}|\nabla_0(R-\overline R)|^2 u^2d \theta +2
\int_0^{2\pi}|R-\overline R|^{2} d \theta + o_t(1)
\int_0^{2\pi}|\nabla_0(R-\overline R)|^{2} d
\theta \nonumber \\
=&(-2+o_t(1))\int_0^{2\pi}|\nabla_0(R-\overline R)|^2 d \theta  +2
\sum_{i=1}^\infty (R^i)^2, \nonumber
\end{align}
where  $o_t(1) \to 0$ as $t \to \infty$, and we use the fact
$u=1+o_t(1)$ in the third identity. On the other hand, from
\begin{equation}
u_{\theta \theta} +\frac 14 u= \frac {R} 4 u^{-3}, \label{3-13}
\end{equation}
we have
$$ (u \cdot u_{\theta \theta})_\theta+\frac
14(u^2)_\theta=\frac 14R_\theta u^{-2} -\frac 12 R u^{-3} \cdot
u_\theta.$$ Using (\ref{3-13}) one more time, we have
$$
 (u \cdot u_{\theta \theta})_\theta+\frac
14(u^2)_\theta+2 u_\theta \cdot u_{\theta \theta} +\frac 12 u
\cdot u_\theta =\frac14 R_\theta u^{-2}.$$ Thus we have a
pointwise  Kazdan-Warner type  identity:
$$
(u^2)_{\theta \theta \theta} +(u^2)_\theta =\frac{ R_\theta}2
u^{-2}.
$$
 Hence for $i=1, 2$,
$$
\int_0^{2\pi} \frac{ R_\theta}2 u^{-2}\phi_i d \theta =
(\mu_i^3-\mu_i) \int_0^{2\pi} u^{-2} \phi_i d\theta=0,
$$
and
\begin{align*}
R^i& = \int_0^{2\pi} (R-\overline{R}) \cdot \phi_i d \theta
=-\frac{1}{\mu_i}\int_0^{2\pi} (R-\overline{R}) \cdot \Delta_0\phi_i d \theta
= \int_0^{2 \pi} R_\theta \cdot (\phi_i)_\theta d \theta \\
&=\int_0^{2\pi}R_\theta\cdot(\phi_i)_\theta u^{-2}d\theta+\int_0^{2\pi}
R_\theta \cdot (\phi_i)_\theta (1-u^{-2})d \theta\\
&=0+\int_0^{2 \pi} R_\theta \cdot(\phi_i)_\theta(1-u^{-2})d\theta
= o_t(1)\int_0^{2 \pi}|\nabla_0(R-\overline{R})|d\theta\\
&= o_t(1)\left(\int_0^{2
\pi}|\nabla_0(R-\overline{R})|^2d\theta\right)^{1/2}.
\end{align*}
It follows from (\ref{3-12}) that
\begin{align*}
\partial_t \int_0^{2\pi}|R-\overline R|^{2} d \sigma&=(-2+o(1))
\int_0^{2\pi}|\nabla_0(R-\overline R)|^2 d \theta\\
&\qquad\qquad  +2 \sum_{i=3}^\infty (R^i)^2 +2((R^1)^2+(R^2)^2)\\
&=(-2+o(1))\int_0^{2\pi}|\nabla_0(R-\overline R)|^2 d\theta+2\sum_{i=3}^\infty(R^i)^2\\
&=(-2+o(1))\sum_{i=1}^\infty \mu_i (R^i)^2+2\sum_{i=3}^\infty(R^i)^2\\
&\le(-2+o(1))\sum_{i=1}^\infty(R^i)^2.
\end{align*}
Therefore, we have, for $t$ sufficiently large,
$$
\partial_t \int_0^{2\pi}|R-\overline R|^{2} d\sigma\le
-a\int_0^{2\pi}|R-\overline R|^{2} d\sigma
$$
for some positive constant $a$. This implies
\begin{equation}\label{3-16}
\int_0^{2\pi}|R-\overline R|^{2} d \sigma\le Ce^{-at} \end{equation}
 for
some positive constant $C$.

For any $T>0$ and $\delta\in[0,1]$, integrating (\ref{3-10})
from $T$ to $T+\delta$ and using (\ref{3-16}) we have
$$
\int_T^{T+\delta}\left(\int_0^{2\pi}|\nabla(R-\overline{R})|^{2}d\sigma\right)dt
\le Ce^{-aT},
$$
which implies
\begin{align*}
\int_T^{T+\delta}||R-\overline{R}||_{L^\infty}dt\le&\sqrt{2\pi}
\int_T^{T+\delta}\left(\int_0^{2\pi}|\nabla(R-\overline{R})|^{2}
d\sigma\right)^{\frac12} dt\\
\le& 2\pi\left(\int_T^{T+\delta}\int_0^{2\pi}|\nabla(R-\overline{R})|^{2}
d\sigma dt \right)^{\frac12}\le C e^{-aT/2}.
\end{align*}
Note that along the Yamabe flow (\ref{Yequ}), $(\ln u)_t=\frac14(R
-\overline{R})$. Integrating from $T$ to $T+\delta$ we obtain that
$$
|\ln u(\theta,T+\delta)-\ln u(\theta, T)|\le Ce^{-aT/2},
$$
for any $\theta\in[0,2\pi]$, $T>0$ and $\delta\in[0,1]$. Hence
$||u(t)-u_\infty||_{L^\infty}\le C e^{-at/2}$. We hereby complete
the proof of Theorem 1.
\medskip

\section{convergence for the affine flow}

We shall focus on the proof of the exponential convergence for the
affine flow in this section. It is not clear to us whether one can
obtain the $L^\infty$ estimate on the metric via the method of
moving plane. Instead, we first obtain the $L^\infty$ estimate on
the curvature via the integral estimates, then derive the
exponential convergence for the curvature along the same line as
in \cite{C1}. The exponential convergence of the metric follows
immediately.

 Along the affine flow (\ref{Aequ}), $u$ satisfies
\begin{equation}\label{4-1}
u_t=\frac14(\kappa-\overline{\kappa }) u,
\mbox{ on }S^1\times[0,\infty),
\end{equation}
and $\kappa_g$ satisfies
\begin{equation}\label{4-2}
\kappa_t=\frac14\Delta\kappa +\kappa (\kappa -\overline{\kappa }),
\end{equation}
where and throughout  this section, we use $\kappa$ to replace
$\kappa_g$ and $\Delta$ to replace $\Delta_g$. As in the Yamabe
flow case, we define: for $p\ge 2$,
\begin{equation}\label{4-4}
\mathcal{F}_p(t)=\int_0^{2\pi} |\kappa -\overline{\kappa }|^p d \sigma.
\end{equation}
Then $\mathcal{F}_p(t)$ satisfies
\begin{align}\label{atfp}
\partial_t \mathcal{F}_p= &-\frac {p-1}{p} \int_0^{2\pi}(|\kappa
 -\overline\kappa |^{p/2})_\sigma^2 d \sigma + (p-\frac 12)\int_0^{2\pi}|\kappa
-\overline \kappa |^{p}(\kappa -\overline \kappa ) d \sigma \\
&+p \overline \kappa  \int_0^{2\pi}|\kappa -\overline \kappa |^{p} d \sigma-
\frac p{4\pi}\int_0^{2\pi}|\kappa -\overline \kappa |^{p-2}  (\kappa
-\overline \kappa ) d \sigma \cdot\int_0^{2\pi}|\kappa
-\overline \kappa |^{2} d \sigma\nonumber\\
&\le C_1(p)(\mathcal{F}_{p+1}+\mathcal{F}_p+\mathcal{F}_p^{1+\frac1p})
-\frac{p-1}p\int_0^{2\pi}|\nabla(|\kappa
 -\overline\kappa |^{p/2})|^2 d \sigma. \nonumber
\end{align}
Similar arguments to that in previous section imply that
\begin{lemma}\label{lem4-1} For any $p \ge 2$,
$$
\mathcal{F}_p(t) \to 0 \ \mbox{ as }\  t \to +\infty \ \ \mbox{ and } \ \
\int_0^\infty \mathcal{F}_p(t)dt<\infty.
$$
\end{lemma}
Next, we shall obtain $L^\infty$ convergence of $\kappa$. The
argument is different to that in previous section, since we do not
know that $u(\theta, t)$ is uniformly bounded below and above by
positive constants.

\begin{lemma}
$$||\kappa-\overline{\kappa}||_{L^\infty}\to 0 \mbox{ as }t\to\infty.$$
\end{lemma}
\begin{proof}
Direct computation yields
\begin{align}
\label{4-5}\partial_t\int_0^{2\pi} (\kappa_\sigma)^2d\sigma
=&-\frac12\int_0^{2\pi}(\kappa_{\sigma\sigma})^2d\sigma
+\frac12\int_0^{2\pi}(\kappa_\sigma)^2(9\kappa-5\overline\kappa)d\sigma\\
\le&
\frac12\int_0^{2\pi}(\kappa_\sigma)^2(9\kappa-5\overline\kappa)
d\sigma.
\nonumber
\end{align}
Set $p=2$ and $p=3$ in (\ref{atfp}), we obtain that
$$
\int_0^\infty\int_0^{2\pi}\kappa_\sigma^2d\sigma dt<\infty\mbox{ and }
\int_0^\infty\int_0^{2\pi}|\kappa-\overline{\kappa}|\kappa_\sigma^2
d\sigma dt<\infty.
$$
Therefore, for any $\epsilon>0$, there exists $t_\epsilon>0$, such that
$$
\int_0^{2\pi}\kappa_\sigma^2(\sigma,t_\epsilon)d\sigma <\epsilon
\mbox{ and }\int_{t_\epsilon}^\infty\int_0^{2\pi}(\kappa_\sigma)^2
(9\kappa-5\overline\kappa)d\sigma dt<\epsilon.
$$
Then for any $t>t_\epsilon$, integrating (\ref{4-5}) from
$t_\epsilon$ to $t$ we obtain that
$$
\int_0^{2\pi}\kappa_\sigma^2(\sigma,t)d\sigma <
\int_0^{2\pi}\kappa_\sigma^2(\sigma,t_\epsilon)d\sigma+\epsilon<2\epsilon.
$$
Hence
$$
\lim_{t\to\infty}\int_0^{2\pi} (\kappa_\sigma)^2d\sigma=0\mbox{
and } \int_0^\infty\int_0^{2\pi}\kappa_{\sigma\sigma}^2d\sigma
dt<\infty.
$$
For any $t>0$, choose $\sigma_0$ s.t. $\kappa(\sigma_0)=\overline{\kappa}$.
Then we have
\begin{equation}\label{infty}
|\kappa(\sigma)-\overline{\kappa}|=\left|\int_{\sigma_0}^\sigma
(\kappa-\overline{\kappa})_\sigma d\sigma\right|\le\int_0^{2\pi}
|\kappa_\sigma| d\sigma\le\left(2\pi\int_0^{2\pi} (\kappa_\sigma)^2d\sigma
\right)^{1/2}\to 0.
\end{equation}
\end{proof}

Because of the correspondence between positive functions $u$
satisfying the orthogonal condition (\ref{on}) and convex simple
closed curves ${\bf C}_u\subset{\bf R}^2$, we may define the
Euclidean perimeter $\ell(g)$ of $(S^1, g)$,  and the Euclidean
area  $A(g)$ of $(S^1, g)$, to be the Euclidean perimeter of ${\bf
C}_u$ and the Euclidean area of ${\bf C}_u$, respectively.

\begin{proposition}
\label{propa} Along the flow (\ref{Aequ}), the area $A(g(t))$ is a
decreasing function.
\end{proposition}
To find the evolution equation for $A(g(t))$, the Euclidean area
of $g(t)$, we consider the support function $h(\theta,t)$ of the
corresponding curves ${\bf x}(\theta,t)$.

In general, the support function  $h: S^1\to \bf R$, of a compact
convex set $C\subset{\bf R}^2$ is defined by
$$
h(\theta)=\max_{x\in C} \{(\cos\theta,\sin\theta)\cdot x\},
$$
where $``\cdot "$ is the usual dot product in ${\bf R}^2$. The
support function of a convex simple closed curve is defined to be
the support function of the set enclosed by the curve. And we
define the support function of $u$ to be the support function of
corresponding curve ${\bf C}_u$. Thus the support function of
${\bf C}_u$ is given by
$$
h(\theta,t)=\sin\theta\int_0^\theta u^{-3}(\phi,t)\cos\phi d\phi
-\cos\theta\int_0^\theta u^{-3}(\phi,t) \sin\phi d\phi,
$$
and satisfies the following evolution equation:
$$
\partial_t h(\theta,t)=\frac34(\overline{\kappa}h(\theta,t)
-u(\theta,t)+u(0,t)\cos\theta +u_\theta(0,t)\sin\theta).
$$

\begin{lemma}
\label{adec}
$A(g(t))$ satisfies the following evolution equation:
\begin{equation}\label{at}
\partial_t A(g(t))=\frac32\overline{\kappa}(t)A(g(t))-\frac32 \pi.
\end{equation}
\end{lemma}
\begin{proof}
Note that the Euclidean curvature of ${\bf C}_u(\theta,t)$ is
$u^3(\theta,t)$. We have
$$
A(g(t))=\frac12\int_0^{2\pi}h(\theta,t)u^{-3}(\theta,t)d\theta.
$$
It follows that
\begin{align*}
\partial_t A(g(t))=&\frac12\int_0^{2\pi}\frac34(\overline{\kappa}
h-u)u^{-3}d\theta+\frac12\int_0^{2\pi}(-\frac34)h u^{-3}(\kappa
-\overline{\kappa})d\theta\\
=&\frac{3}{4}\overline{\kappa}(t)\int_0^{2\pi}hu^{-3}d\theta
-\frac38\int_0^{2\pi}u^{-2}d\theta-\frac38\int_0^{2\pi}h(u_{\theta\theta}
+u)d\theta\\
=&\frac{3}{2}\overline{\kappa}(t)A(g(t))-\frac{3}{4}\pi
-\frac38\int_0^{2\pi}u(h_{\theta\theta}+h)d\theta\\
=&\frac{3}{2}\overline{\kappa}(t)A(g(t))-\frac32\pi,
\end{align*}
where we used that $\int_0^{2\pi}u^{-2}d\theta=2\pi$ and
$h_{\theta\theta}+h =u^{-3}$.
\end{proof}

Note that $\overline{\kappa}$ is a  strictly increasing function
of $t$. If at some $t_0$, $\overline{\kappa}(t_0)A(g(t_0))\ge
\pi$, (\ref{at}) immediately implies that $A(g(t))$ increases for
$t\ge t_0$ and $A(g(t))\to\infty$ as $t\to\infty$. On the other
hand, it is known (see, e.g. \cite{su} Theorem I, page 48)
$$
\pi(\min\kappa)^{-3/2}\le A(g(t))\le\pi(\max\kappa)^{-3/2},
$$
thus $A(g(t))$ is bounded as $t\to\infty$. Contradiction.
Therefore
$$
\partial_t A(g(t))=\frac32\overline{\kappa}(t)A(g(t))-\frac32 \pi<0,
\mbox{ for all }t\ge 0.
$$
This completes the proof of Proposition \ref{propa}.

\smallskip

For any $t\ge 0$, as in Section 3 of \cite{andrews}, we may apply
optimal special linear transformations to ${\bf C}_{u(t)}$ so that
its  arc length is minimized. We denote the resultant modified
curve solution by $\tilde{\bf C}_{u(t)}$ and the corresponding
modified solution by $\tilde{g}(t)$.
\begin{proposition}
\label{propg} The modified solution $\tilde{g}(t)\to g_s$ in
$L^\infty (S^1)$ as $t\to\infty$.
\end{proposition}

For a positive function $u(\theta)$ defined on $S^1$, $\lambda>0$
and $\alpha\in[0,2\pi]$, we  consider transformations
$(T_{\lambda,\alpha}
u)(\theta):=u(\sigma_{\lambda,\alpha}(\theta))
\psi_{\lambda,\alpha}(\theta)$, where
$$
\psi_{\lambda,\alpha}(\theta)=\sqrt{\lambda\cos^2(\theta-\alpha)
+\lambda^{-1}\sin^2(\theta-\alpha)},
$$
and
$$
\sigma_{\lambda,\alpha}(\theta)=\alpha+\int_0^\theta
\psi_{\lambda,\alpha}^{-2}d\theta=\alpha+ \left\{
\begin{array}{ll}
\arctan(\lambda^{-1}\tan(\theta-\alpha)), & \theta-\alpha \in [0, \frac\pi2]\\
\arctan(\lambda^{-1}\tan(\theta-\alpha))+\pi,
& \theta-\alpha \in (\frac\pi2, \frac{3\pi}2]\\
\arctan(\lambda^{-1}\tan(\theta-\alpha))+2 \pi,
& \theta-\alpha \in (\frac{3\pi}2, 2 \pi).
\end{array}
\right.
$$
If $u$ satisfies the orthogonal condition (\ref{on}), the
corresponding curve ${\bf C}_u(\theta)\subset {\bf R}^2$ is a
convex simple closed curve. Then the transformations $\tilde{T}$
on ${\bf C}_{u}$ corresponding to $T$ are special linear
transformations. To modify the curve ${\bf C}_u$ so that its
length is minimized, is equivalent to find $\lambda,\alpha$ so
that $\ell(T_{\lambda_u,\alpha_u}u)$ is minimized. Since as
$\lambda\to 0$ or $\lambda\to\infty$,
$\ell(T_{\lambda,\alpha}u)\to\infty$, we can  always find a pair
($\lambda_u, \alpha_u$) so that $\ell(T_{\lambda_u,\alpha_u}u)$ is
minimized. For such a pair, we have

\begin{lemma}
\begin{equation}\label{critical}
\int_0^{2\pi}\frac{\cos2\theta}{(T_{\lambda_u,\alpha_u}u)^3}d\theta
=\int_0^{2\pi}\frac{\sin2\theta}{(T_{\lambda_u,\alpha_u}u)^3}d\theta=0.
\end{equation}
\end{lemma}
\begin{proof}
\begin{align*}
\ell(T_{\lambda,\alpha}u)=&\int_0^{2\pi}(T_{\lambda,\alpha}u)^{-3}d\theta
=\int_0^{2\pi}\frac{1}{u^3(\sigma_{\lambda,\alpha})(\theta)
\psi^3_{\lambda,\alpha}(\theta)}d\theta\\
=&\int_0^{2\pi}\frac{\sqrt{\lambda\sin^2(\phi-\alpha)
+\lambda^{-1}\cos^2(\phi-\alpha)}}{u^3(\phi)}d\phi.
\end{align*}
At $(\lambda_u,\alpha_u)$, we have
$$
\frac{\partial}{\partial \lambda}\ell(T_{\lambda,\alpha}u)=0\mbox{ and }
\frac{\partial}{\partial \alpha}\ell(T_{\lambda,\alpha}u)=0,
$$
which, by direct computation, imply (\ref{critical}).
\end{proof}
Write $v=T_{\lambda_{u},\alpha_{u}}u$. Then the modified curve
$\tilde{\bf C}_u={\bf C}_v$.
\begin{lemma}
\label{lbdd}
$$
\int_0^{2\pi}v^{-3}(\theta)d\theta\le 2\sqrt{6\pi A(v)},
$$
where $A(v)$ is the area of the region enclosed by ${\bf C}_v$.
\end{lemma}
\begin{proof}
Suppose the Fourier expansion of the support function $h_v$ of $v$
is $h_v(\theta)=a_0+\sum_{n=1}^{\infty}a_n\cos
n(\theta-\theta_n)$.
 Then (\ref{critical})
implies that $a_2=0$, thus
$$
v^{-3}=h_v+(h_v)_{\theta\theta}=a_0+\sum_{n=3}^{\infty}a_n(1-n^2)
\cos n(\theta-\theta_n).
$$
By integrating nonnegative function $v^{-3}(\theta)(1\pm\cos
n(\theta-\theta_n))$ on $[0,2\pi]$,  we obtain that $2\pi a_0\pm
\pi a_n(1-n^2)\ge 0$. It follows that $|a_n|\le 2a_0/(n^2-1)$ for
all $n\ge3$.
\begin{align*}
A(v)=&\frac12\int_0^{2\pi}h_v(\theta) v^{-3}(\theta)d\theta
=\pi a_0^2+\frac\pi2\sum_{n=3}^{\infty}(1-n^2)a_n^2\\
\ge&\pi
a_0^2\left(1-\frac12\sum_{n=3}^{\infty}\frac4{n^2-1}\right)
=\frac16\pi a_0^2.
\end{align*}
Therefore $a_0\le \sqrt{6 A(v)/\pi}$  and
$$
\int_0^{2\pi}v^{-3}d\theta=\int_0^{2\pi}(h_v+(h_v)_{\theta\theta})
d\theta=2\pi a_0\le 2\sqrt{6\pi A(v)}.
$$
\end{proof}

 Suppose that $g(t)=u^{-4}(\theta,t)g_s$ is a
solution to the affine flow equation (\ref{affine}) and ${\bf
C}_{u(t)}$ is the corresponding curve solution. For each $t>0$,
write $v(t)=T_{\lambda_{u(t)},\alpha_{u(t)}}u(t)$. Then the
modified curve solution $\tilde{\bf C}_{u(t)}={\bf C}_{v(t)}$ and
$\tilde{g}=v^{-4}g_s$. Since special linear transformations does
not chang the area, we have $A(g(t))=A(\tilde{g}(t))$. Denote the
1-scalar curvature of $\tilde{g}$ by $\tilde{\kappa}$. We have
\begin{equation}\label{rin}
\tilde{\kappa}(\theta,t)=v(\theta,t)^3(v_{\theta\theta}(\theta,t)
+v(\theta,t))=\kappa(\sigma_{\lambda_{u(t)},\alpha_{u(t)}}(\theta),t).
\end{equation}
The following lemma immediately implies Proposition \ref{propg}.
\begin{lemma}
\label{lemma11}
$$
\lim_{t\to\infty}v(\theta,t)=1.
$$
\end{lemma}
\begin{proof}
It follows from Lemma \ref{adec} and Lemma \ref{lbdd} that
$$
\int_0^{2\pi}v^{-3}(\theta,t)d\theta\le 2\sqrt{6\pi A(g(t))}\le
2\sqrt{6\pi A(g(0))}.
$$
It also follows from the fact that $\kappa$ is bounded, thus
$A(t)$ is bounded, which yields the boundedness for
$\int_0^{2\pi}v^{-3}(\theta,t)d\theta$.  Since
$\overline{\kappa}(t)$ is increasing and bounded above by $1$, we
may assume that
$\lim_{t\to\infty}\overline{\kappa}(t)=\kappa_\infty$. It follows
from (\ref{rin}) that
$\lim_{t\to\infty}\tilde{\kappa}(\theta,t)=\kappa_\infty$.
Therefore
$$
\int_0^{2\pi}v^{-3}(\theta,t)\tilde{\kappa}(\theta,t)(\theta,t)d\theta
=\int_0^{2\pi}v(\theta,t)d\theta
$$
is bounded by some constant $c=c(g(0))$. Furthermore, for any
$\alpha, \beta$
$$
\int_\alpha^{\beta}v^{-3}(\theta,t)\tilde{\kappa}(\theta,t)(\theta,t)d\theta
=\int_\alpha^\beta
v(\theta,t)d\theta+v_\theta(\beta,t)-v_\theta(\alpha,t).
$$
By choosing $\alpha$ being the critical point of $v(\theta, t)$
for any fixed $t$, we know from the above  that
$v_\theta(\theta,t)$ is bounded by some constant $c(g(0))$.
Therefore $v(\theta,t)$ is uniformly bounded in $H^1(S^1)$ for all
$t\ge 0$. Hence there exist a sequence $t_n\to\infty$, such that
$v(\theta,t_n)\rightharpoonup v_\infty(\theta)$ in $H^1(S^1)$.
From Sobolev embedding theorem we have $v(\theta,t_n) \to
v_\infty(\theta)$ in $C^{0,\alpha}$ for any
$\alpha\in(0,\frac12)$. Since
$\int_0^{2\pi}v^{-3}(\theta,t)d\theta$ is bounded, we obtain that
$v_\infty(\theta)>0$ and $v_\infty$ satisfies
$v_\infty^3((v_\infty)_{\theta\theta}+v_\infty)=\kappa_\infty$.
Note that $\int_0^{2\pi} v^{-2}(\theta,t)d\theta =\int_0^{2\pi}
u^{-2}(\theta,t)d\theta=1$, it follows from  Theorem 2 in
\cite{NZ1} (Theorem B in the appendix) that $\kappa_\infty=1$. Due
to the orthogonal condition (\ref{critical}) we know  that
$v_\infty=1$. Using the same argument we can prove that any
convergent subsequence of $v(\theta,t)$ converges to $1$. Since
Hence $v(\theta,t)$ is uniformly bounded in $H^1(S^1)$, we have
$\lim_{t\to \infty}v(\theta,t)=1$.
\end{proof}
Now we are ready to prove the exponential convergence. For any
given $t$, define variable $\gamma$ as the inverse of $\theta$
under map $\sigma_{\lambda_{u(t)}, \alpha_{u(t)}}$, that is
$\sigma_{\lambda_{u(t)},\alpha_{u(t)}}(\gamma)=\theta$. Then we
have
$$
d\sigma=u(\theta)^{-2}d\theta=v(\gamma)^{-2}d\gamma.$$ Since the
curvature of metric $d \gamma \otimes d \gamma$ is 1, we know $$
\kappa(\theta(\gamma))=v(\gamma)^3(v''(\gamma)+v(\gamma)).
$$
Suppose that
$$
\kappa=\overline{\kappa}+\sum_{n=1}^\infty(a_n\cos(n\sigma)+b_n\sin(n\sigma))
=\tilde{c}+\sum_{n=1}^\infty(\tilde{a}_n\cos(n\gamma)+\tilde{b}_n\sin(n\gamma)).
$$
Since $v(t)\to 1$ as $t\to\infty$(Lemma \ref{lemma11}), we obtain that
\begin{equation}\label{ab}
a_n=\tilde{a}_n+o(1){\mathcal F}^{\frac12}_2, \quad
b_n=\tilde{b}_n+o(1){\mathcal F}^{\frac12}_2,\qquad
n=0,1,2,3,\cdots,
\end{equation}
where $a_0=\overline {\kappa}$ and $\tilde {a}_0=\tilde c.$ And
\begin{align*}
\tilde{a}_1=&\frac1\pi\int_0^{2\pi}\kappa\cos\gamma d\gamma
=\frac1\pi\int_0^{2\pi}v^3(\gamma)(v''(\gamma)+v(\gamma))\cos\gamma d\gamma\\
=&\frac1\pi\int_0^{2\pi}(v^3(\gamma)-1)(v''(\gamma)+v(\gamma))\cos\gamma
d\gamma\\
=&\frac1\pi\int_0^{2\pi}(1-v^{-3}(\gamma))\kappa \cos\gamma
d\gamma\\
=&\frac1\pi\int_0^{2\pi}(1-v^{-3}(\gamma))(\kappa-
\overline{\kappa}) \cos\gamma d\gamma=o(1){\mathcal
F}^{\frac12}_2,
\end{align*}
where we use the orthogonality of $v$. Similarly
$\tilde{b}_1=o(1){\mathcal F}^{\frac12}_2$.

To estimate $\tilde{a}_2$ and $\tilde{b}_2$, we note that
$\kappa=v^3(\gamma)(v''(\gamma)+v(\gamma))$ implies
$2\kappa_\gamma v^{-2}=(v^2)_{\gamma\gamma}+4(v^2)_\gamma$. Hence
$$
\int_0^{2\pi}\kappa_\gamma\sin(2\gamma) v^{-2} d\gamma
=\int_0^{2\pi}\kappa_\gamma\cos(2\gamma) v^{-2} d\gamma=0.
$$
It follows that
\begin{align*}
\tilde{a}_2&=\frac1\pi\int_0^{2\pi}\kappa\cos(2\gamma) d\gamma
=-\frac1{2\pi}\int_0^{2\pi}\kappa_\gamma\sin(2\gamma) d\gamma\\
=&-\frac1{2\pi}\int_0^{2\pi}\kappa_\gamma\sin(2\gamma) v^{-2} d\gamma
-\frac1{2\pi}\int_0^{2\pi}\kappa_\gamma\sin(2\gamma)(1- v^{-2}) d\gamma\\
=&o(1)\left(\int_0^{2\pi}\kappa_\gamma^2d\gamma\right)^{\frac12}
=o(1)\left(\int_0^{2\pi}\kappa_\sigma^2d\sigma\right)^{\frac12}.
\end{align*}
Similarly $\tilde{b}_2=o(1)\left(\int_0^{2\pi}\kappa_\sigma^2
d\sigma\right)^{\frac12}$. It follows from (\ref{ab}) that
$$
a_1,a_2,b_1,b_2=o(1)\left(\int_0^{2\pi}\kappa_\sigma^2d\sigma\right)^{\frac12}.
$$
Therefore, from (\ref{atfp}) we have
\begin{align*}
\partial_t{\mathcal F}_2=&-\frac12\int_0^{2\pi}\kappa_\sigma^2d\sigma
+\frac32\int_0^{2\pi}(\kappa-\overline{\kappa})^3d\sigma+2\overline{\kappa}
\int_0^{2\pi}(\kappa-\overline{\kappa})^2d\sigma\\
=&(-\frac12+o(1))\int_0^{2\pi}\kappa_\sigma^2d\sigma+2
\int_0^{2\pi}(\kappa-\overline{\kappa})^2d\sigma\\
=&(-\frac12+o(1))\int_0^{2\pi}\kappa_\sigma^2d\sigma
+2\sum_{n=1}^2(a_n^2+b_n^2)+2\sum_{n=3}^{\infty}(a_n^2+b_n^2)\\
\le&(-\frac12+o(1)){\mathcal F}_2.
\end{align*}
This implies
\begin{equation}\label{kappaexp}
{\mathcal F}_2=
\int_0^{2\pi}(\kappa-\overline{\kappa})^{2}d\sigma\le Ce^{-at},
\end{equation}
for some positive constant $a$ and $C$.

Using the same argument as in the proof of Theorem 1, we obtain that
$$
|\ln u(\theta,T+\delta)-\ln u(\theta, T)|\le Ce^{-aT/2},
$$
for any $\theta\in[0,2\pi]$, $T>0$ and $\delta\in[0,1]$. Hence
$\lim_{t\to\infty}u(\theta,t)=u_\infty(\theta)$, with
$||u(t)-u_\infty||_{L^\infty}\le C e^{-at/2}$ and the 1-curvature
of $g_\infty:=u_\infty^{-4}g_s$ is constant $1$. This completes
the proof of Theorem \ref{theorem2}.

\section{appendix: Sobolev and  Blaschke-Santal\'o inequality}
For readers' convenience, we list here the basic theorems on sharp
embedding inequalities used in the paper. The proofs of these
theorems can be found in, e.g. \cite{NZ1}.

\noindent {\bf Theorem A} (General Blaschke-Santal\'o): For
$u(\theta)\in H^1(S^1)$ and $u>0$, if $u$ satisfies
$$
\int_0^{2\pi}\frac{\cos(\theta+\gamma)}{u^3(\theta)}d\theta =0
$$
for all $\gamma$, then
$$
\int_0^{2\pi}(u^2_\theta-u^2)d\theta
\int_0^{2\pi}u^{-2}(\theta)d\theta\ge -4\pi^2,
$$
and the equality holds if and only if
$$
u(\theta)=c\sqrt{\lambda^2\cos^2(\theta-\alpha)+\lambda^{-2}\sin^2(\theta-\alpha)},
$$
for some $\lambda,c>0$ and $\alpha \in [0, 2\pi)$. \label{thmA-1}

\medskip

\noindent{\bf Theorem B}: For $u(\theta)\in H^1(S^1)$ and $u>0$,
$$
\int_0^{2\pi}(u^2_\theta-\frac14u^2)d\theta
\int_0^{2\pi}u^{-2}(\theta)d\theta\ge -\pi^2,
$$
and the equality holds if and only if
$$
u(\theta)=c\sqrt{\lambda^2\cos^2\frac{\theta-\alpha}2
+\lambda^{-2}\sin^2\frac{\theta-\alpha}2}
$$
for some $\lambda,c>0$  and $\alpha \in [0, 2\pi)$. \label{thmA-2}

\medskip
\noindent{\bf ACKNOWLEDGMENT. } The work of M.Zhu is partially
supported by the NSF grant DMS-0604169.


\begin{thebibliography}{99}
\bibitem{Al}Alvarez, Luis; Guichard, Fr¨¦d¨¦ric; Lions, Pierre-Louis;
Morel, Jean-Michel, "Axioms and fundamental equations of image
processing". {\it Arch. Rational Mech. Anal.} 123 (1993), no. 3,
199--257.
\bibitem{A1} Andrews, B., Contraction of convex hypersurfaces by their affine normal.
 J. Differential Geom. 43 (1996), no. 2, 207--230.
\bibitem{andrews}Andrews, B., The Affine Curve-Lengthening Flow, J. Reine Angew.
Math. 506 (1999), 43-83.
\bibitem{A2}Andrews, B., Classification of limiting shapes for
isotropic curve flows. J. Amer. Math. Soc. 16 (2003), no. 2,
443--459.
\bibitem{CGS}Caffarelli, L., Gidas, B.  and  Spruck, J.,
Asymptotic symmetry and local behavior of semilinear elliptic
equations with critical Sobolev growth, Comm. Pure Appl. Math.
 42 (1989),  271-297.
\bibitem{C1}Chen, X.X., Calabi Flow in Riemann Surfaces Revisited: A New Point
of View, Internat. Math. Res. Notices, 2001, no. 6, 275-297.
\bibitem{G}Gage, M., An isoperimetric inequality
with applications to curve shortening. Duke Math. J. 50 (1983), no. 4, 1225--1229.
\bibitem{GH} Gage, M.; Hamilton, R. S., The heat equation shrinking
convex plane curves. J. Differential Geom. 23 (1986), no. 1, 69--96.
\bibitem {NZ1}Ni, Y. and Zhu, M., Steady States for One Dimensional Conformal
Metric Flows, \url{http://arxiv.org/abs/math.AP/0611254}
 \bibitem{Si} Simon, L.,  Asymptotics for a class of nonlinear evolution equations,
  with applications to geometric problems. Ann. of Math. (2) 118 (1983), no. 3, 525--571
\bibitem{S}Struwe,M.,Curvature flows on surfaces. Ann. Sc. Norm. Super.
Pisa Cl. Sci. (5) 1 (2002), no. 2, 247--274
\bibitem{ST}Sapiro, G. and Tannenbaum, A., On Affine Plane Curve Evolution,
J. Funct. Anal., 119(1994), 79-120.
\bibitem{su}Su, B, Affine differential geometry, Gordon and Breach,
New York, 1983.
\bibitem {Ye}Ye, R., Global Existence and Convergence of Yamabe Flow,
J. Differential Geom., 39(1994), 35-50.
\end{thebibliography}
\end{document}